\UseRawInputEncoding

\documentclass[11pt,letterpaper]{amsart}
\usepackage{amssymb}
\usepackage[all]{xy}
\usepackage[toc,page,title,titletoc,header]{appendix}

\usepackage{xcolor}

\theoremstyle{plain}
\newtheorem{thm}{Theorem}[section]

\newtheorem*{thm1.1}{Theorem 1.1}

\newtheorem{lem}[thm]{Lemma}
\newtheorem{cor}[thm]{Corollary}
\newtheorem{pro}[thm]{Proposition}

\theoremstyle{definition}

\newtheorem{rem}[thm]{Remark}

\newtheorem{defi}[thm]{Definition}
\newtheorem{exe}[thm]{Example}

\newtheorem{fact}[thm]{Fact}

\numberwithin{equation}{section}
\newcounter{elno}                

\newcommand{\la}{\lambda}

\newcommand{\Supp}{{\rm Supp}}

\newcommand{\Psef}{{\rm Psef}}
\newcommand{\CH}{{\rm CH}}

\newcommand{\Char}{{\rm char}}
\newcommand{\Tr}{{\rm Tr}}

\newcommand{\id}{{\rm id}}

\newcommand{\Gal}{{\rm Gal}}

\newcommand{\Sp}{{\rm Sp}\,}

\newcommand{\boxtensor}{{\Box\kern-9.03pt\raise1.42pt\hbox{$\times$}}}

\newcommand{\propsubset}
{\mbox{$\textstyle{
			\subseteq_{\kern-5pt\raise-1pt\hbox{\mbox{\tiny{$/$}}}}}$}}

\newcommand{\sC}{{\mathcal C}}

\newcommand{\C}{{\mathbb C}}

\newcommand{\F}{{\mathbb F}}

\renewcommand{\P}{{\mathbb P}}
\newcommand{\Q}{{\mathbb Q}}
\newcommand{\R}{{\mathbb R}}

\newcommand{\Z}{{\mathbb Z}}
\newcommand{\bk}{{\mathbf{k}}}

\newcommand{\Fix}{\mathrm{Fix}}

\begin{document}
	\title[]{Numerical spectrums control Cohomological spectrums}
	\author{Junyi Xie}
	
	\address{Beijing International Center for Mathematical Research, Peking University, Beijing 100871, China}
	
	\email{xiejunyi@bicmr.pku.edu.cn}

	\thanks{The author is supported by the NSFC Grant No.12271007.}

	\date{\today}

	\bibliographystyle{alpha}

	
	\maketitle
	

	\begin{abstract}
		Let $X$ be a smooth irreducible projective variety over a field $\bk$ with $\dim X=d.$ Let $\tau: \Q_l\to \C$ be any field embedding.
		Let $f: X\to X$ be a surjective endomorphism. We show that for every $i=0,\dots,2d$, the spectral radius of $f^*$ on the numerical group $N^i(X)\otimes \R$ and on the $l$-adic cohomology group $H^{2i}(X_{\overline{\bk}},\Q_l)\otimes \C$ are the same. As a consequence, if $f$ is $q$-polarized for some $q>1$, we show that the norm of every eigenvalue of $f^*$ on the $j$-th cohomology group is $q^{j/2}$ for all $j=0,\dots, 2d.$ This generalizes Deligne's theorem for Weil's  Riemann Hypothesis to arbitary polarized endomorphisms and proves a conjecture of Tate. 
		We also get some applications for the counting of fixed points and its ``moving target" variant. 
		
		Indeed we studied the more general actions of certain cohomological coorespondences and we get the above results as consequences in the endomorphism setting.
		\end{abstract}
	
	\section{Introduction}
	The aim of the paper is to compare the eigenvalues on numerical groups and on $l$-adic cohomology groups for 
	actions of endomorphisms or more generally for certain cohomological coorespondences on smooth projective varieties.

\subsection{Action of endomorphisms}
Let $X$ be a smooth irreducible projective variety over a field $\bk$ with $\dim X=d.$
Let $f: X\to X$ be a surjective endomorphism.

\subsubsection{Numerical spectrum}
For every $i=0,\dots,d$, the \emph{$i$-th numerical spectrum} of $f$ is 
$\Sp^{num}_i(f)$ the multi-set of eigenvalues of $$f^*: N^i(X)\otimes \R\to N^i(X)\otimes \R.$$
where $N^i(X)$ the group of numerical cycles of codimension $i$ of $X$.
The \emph{$i$-th numerical spectral radius} is  
$$\beta_i(f):=\max \{|a||\,\, a\in \Sp^{num}_i(f)\}.$$
By \cite[Theorem 1.1(3)]{Truong2020} (see also \cite{Dinh2005,Dang2020}), the sequence $\beta_i(f), i=0,\dots, d$ is log-concave i.e. 
$$\beta_i(f)^2\geq \beta_{i-1}(f)\beta_{i+1}(f)$$ for $i=1,\dots, d-1.$

\begin{rem}There is an important notion of \emph{dynamical degrees} $\la_i(f), i=0,\dots, d$, which are well defined even when $f$ is merely a dominant rational self-map. \cite{Russakovskii1997,Dinh2005,Dinh2004,Truong2020,Dang2020}.
In our setting where $f$ is an endomorphism, these two notions coincide i.e. $\la_i(f)=\beta_i(f)$.
\end{rem}
By projection formula and the Poincar\'e duality, we have [c.f. Proposition \ref{prospinf}]
$$\beta_i^-(f):=\min\{|a||\,\, a\in \Sp^{num}_i(f)\}=\deg f/\beta_{d-i}(f).$$

\subsubsection{Cohomological spectrum}
Let $l$ be a prime number  with $l\neq \Char\, \bk.$ 
Let $$G(X):=\{j=0,\dots, 2d|\,\, H^j(X_{\overline{\bk}},\Q_l)\neq 0\}.$$ By Poincar\'e duality, $j\in G(X)$ if and only if $2d-j\in G(X).$
Moreover $$\{2i|\,\, i=0,\dots, d\}\subseteq G(X).$$

\medskip

For $j=0,\dots, 2d$, the $j$-th cohomological spectrum
$\Sp_j(f)$
is the multi-set of eigenvalues of $$f^*: H^j(X_{\overline{\bk}},\Q_l)\to H^j(X_{\overline{\bk}},\Q_l).$$
By Fact \ref{factrationa}, which is a consequence of Deligne's proof of Weil conjecure \cite[Theorem 2]{Katz1974}, all elements in 
$\Sp_j(f)$ are algebraic integers and moreover $\Sp_j(f)$ is $\Gal(\bar{\Q}/\Q)$-invariant.

Fix any embedding $\tau: \Q_l\hookrightarrow \C$ and view $\Q_l$ as a subfield of $\C$ via $\tau.$
The image of $\Sp_j(f)$ in $\C$ does not depend on the choice of $\tau.$ We may view $\Sp_j(f)$ as a subset of $\C$.
For $j=0,\dots, 2d$, the $j$-th \emph{cohomological spectral radius} is 
$$\alpha_j(f):= \max\{|b||\,\, b\in \Sp_j(f)\}.$$ 
In particular, if $j\not\in G(X)$, then $\alpha_j(f)=0.$
For $j\in G(X)$, we have 
$$\alpha^-_j(f):=\min\{|b||\,\, b\in \Sp_j(f)\}=\frac{\deg(f)}{\alpha_{2d-j}(f)}.$$
\subsubsection{Main result and consequences}
The following is our main result in the endomorphism case, which 
is indeed a consequence of the more general result Theorem \ref{thmcpequnpintro} for 
\emph{bi-finite} correspondences.

\begin{thm}[=Corollary \ref{corendosum}]\label{thmendosum}
	For every $i=0,\dots, d$, we have 
	$$\log \alpha^-_{2i}(f)=\log\beta^-_i(f)
	\text{ and }
	\log\beta_i(f)=\log \alpha_{2i}(f).$$
		For every odd $j\in G(X)$, we have 
	$$\frac{\log \beta^-_{j-1}(f)+\log \beta^-_{j+1}(f)}{2}\leq \log \alpha^-_j(f)\leq \log\alpha_j(f)\leq \frac{\log \beta_{j-1}(f)+\log \beta_{j+1}(f)}{2}.$$
	\end{thm}
\begin{rem}
	Theorem \ref{thmendosum} can be state in a more geometric way as follows: 
	Consider the following sets of points in $\R^2$: 
		$$F_{num}:=\{(2i,\log|u|)|\,\,u\in \Sp_i^{num}(f)\} \text{ and } 
	F_{coh}:=\{(j, \log|u|)|\,\, u\in \Sp_j(f))\}.$$
	Let $C_{num}$ and $C_{coh}$ be their convex envelops. Then we have $$C_{num}=C_{coh}.$$
	Moreover for every $i=0,\dots,d$, the end points of $C_{num}\cap (\{2i\}\times \R)$ are contained in $F_{num}.$
	\end{rem}

Theorem \ref{thmcpequnpintro} answers a question of Truong \cite[Question 2]{Truong2016} positively for endomorphisms.
By a tensoring product trick originally due to Dinh \cite[Proposition 5.8]{Dinh2005c}, positive answer of \cite[Question 2]{Truong2016} implies positive answer of \cite[Question 4]{Truong2016}.

\medskip
\subsubsection{Historical notes}
When $\bk=\C$, by comparison theorem between singular and $\Q_l$-cohomologies, we can consider the singular cohomology instead of the $\Q_l$-cohomology. After such translation, Theorem \ref{thmendosum} is well known in complex dynamics (see for example \cite[Section 4]{Dinh2017}). 

In positive characteristic, the first progress was due to Esnault-Srinivas \cite{Esnault2013}. They showed that when $f$ is an automorphism on a smooth projective surface, then $\alpha_2=\beta_1.$
Later Truong \cite{Truong2016} generated Esnault-Srinivas's result in any dimension by showing that 
$$\max_{j=0,\dots, 2d}\alpha_j=\max_{i=0\dots,d}\beta_i.$$ Shuddhodan \cite{Shuddhodan2019} gave an alternative approach towards this equality using dynamical zeta functions. Recently, Hu proved Theorem \ref{thmendosum} for 
endomorphisms of abelian varieties \cite{Hu2019,Hu2019a}.

\medskip

Next we apply Theorem \ref{thmendosum} to endomorphisms satisfying some numerical conditions. 

\begin{defi}
	For $q\geq 1$, 
	we say that $f$ is \emph{$q$-straight}, if $\beta_i(f)=q^i$ for all $i=0,\dots, d.$
	We say that $f$ is \emph{straight} if it is $q$-straight for some $q\geq 1.$
\end{defi}
Recall that an endomorphism $f: X\to X$ is called \emph{$q$-polarized} for some $q>1$, if there is an ample line bundle $L$ such that $f^*L=qL.$
For example, if $k=\F_q$, the $q$-Frobenius is $q$-polarized. 
If $f$ is $q$-polarized, then $\beta_i(f)=q^i$ for $i=0,\dots, d.$ 
So all polarized endomorphisms are straight.
By Theorem \ref{thmendosum}, we get the following consequence.
\begin{cor}\label{corpolarweil}
	Assume that $f: X\to X$ is a $q$-straight endomorphism for some $q\geq 1$.
	Then for every $j\in G(X)$, we have $$\alpha^-_j(f)=\alpha_j(f)=q^{j/2}.$$
\end{cor}

This result generalizes Deligne's theorem \cite{Deligne1974}, which proves Weil's Riemann Hypothesis, for any straight endomorphisms.
However, Deligne's theorem \cite{Deligne1974} plays a key role in our proof of Theorem \ref{thmendosum} (hence Corollary \ref{corpolarweil}). 

\medskip

Our Corollary \ref{corpolarweil} proves a conjecture proposed by Tate in 1964 \cite[\S3, Conjecture (d)]{Tate1964}
in Lecture notes prepared in connection with seminars held at the Summer Institute on Algebraic Geometry, Whitney Estate, Woods Hole, MA, July 6--July 31, 1964.
See also \cite[\S 3, Conjecture (d)]{Tate1965} for the published version. In \cite[Conjecture 1.4]{Hu2021}, Hu and Truong proposed the same conjecture and call it the 
Generalized Weil's Riemann Hypothesis.
When $\bk=\C$ and $f$ is polarized, our theorem was proved by Serre \cite{serre} using the Hodge structure. Serre viewed his result as an K\"ahlerian analogy of Weil's conjecture.

\medskip

An endomorphism $f: X\to X$ is called \emph{int-amplified} if there is an ample line bundle $L$ of $X$ such that $f^*L-L$ is ample \cite{Meng2020}. The following fact was observed by Matsuzawa.
\begin{fact}\cite[Proposition 3.7]{Meng2023b} A surjective endomorphism $f:X\to X$ is int-amplified if and only if $\beta_d(f)>\beta_{d-1}(f).$
\end{fact}
Indeed the log-concavity of $\beta_i(f), i=0,\dots, d$ implies that the sequence $\beta_i(f), i=0,\dots, d$ is strictly increasing in the  int-amplified case.
\medskip

By \cite[Theorem 5.1]{fa}, if $f$ is int-amplified, then for every $n\geq 1$, the set of $n$-periodic points $\Fix(f^n)$ is isolated.
Denote by $\#\Fix(f)$ the number of fixed points counting with multiplicity. Combine the Lefschetz fixed point theorem with Corollary \ref{corpolarweil}, we get the following estimates. 
\begin{cor}\label{corfixpoint}
	Assume that $f: X\to X$ is an int-amplified endmorphism. Then we have 
	$$\#\Fix(f^n)=\beta_d^n+O((\beta_d\beta_{d-1})^{n/2}).$$
	In particular, if $f$ is $q$-straight for some $q>1$, then we have 
	$$\#\Fix(f^n)=q^{dn}+O(q^{(d-1/2)n}).$$
\end{cor}
\begin{rem} When $\bk=\C$, stronger results can be proved. For example, in the recent work \cite{Dinh2023}, Dinh-Zhong can counted the periodic points without multiplicity when $f$ is merely a rational map with dominant topological degree.
\end{rem}

We also proved a  ``moving target" version of  Corollary \ref{corfixpoint}, see Proposition \ref{promovingtaget}.

\begin{pro}[=Proposition \ref{promovingtaget}]\label{promovingtagetintro}Assume that $f: X\to X$ is an int-amplified endmorphism. Let $L$ be an ample line bundle on $X$.
	Let $h_n: X\to X, n\geq 0$ be a sequence of endomorphisms of $X$ with 
	$$\limsup_{n\to \infty}(h_n^*L\cdot L^{d-1})^{1/n}<\beta_{d}(f)/\beta_{d-1}(f).$$
	Then for $n\gg 0$, $f^n$ intersects $h_n$ properly in $X\times X$. Moreover, for every $\epsilon>0$ we have 
	$$\#\{f^n(x)=h(x)\}=q^{dn}+o((\beta_d\beta_{d-1})^{(1+\epsilon)n})$$
	counting with multiplicity.
\end{pro}

\subsection{Actions of cohomological correspondences}
Let $X$ be a smooth irreducible projective variety over a field $\bk$ with pure dimension $\dim X=d.$ Here $X$ may not be irreducible. 

\medskip

Denote by $\sC(X,X)_{\Z}$ the space of degree $0$ cohomological correspondences (for the $\Q_l$-cohomology) from $X$ to itself i.e. the image of the cycle class map  $$cl:\CH^d(X\times X)\to H^d(X_{\overline{\bk}}\times X_{\overline{\bk}})(d).$$  With the composition $\circ$, it forms a ring. Denote by $\sC(X,X):=\sC(X,X)_{\Z}\otimes \Q.$ 

\medskip

For $c\in \sC(X,X)$, we define the numerical and cohomological spectral radius $\beta_i(c), i=0,\dots, d$, $\alpha_j(c), j=0,\dots, d$ in the same way as the case for endomorphisms. See Section \ref{subsectioncohr} and \ref{subsectionnumr} for details.

\medskip

Every endomorphism $f: X\to X$ can be viewed as a $d$-cycle in $X\times X$ via its graph, hence induces a cohomological correspondences.
\begin{defi}A \emph{finite correspondence} (resp. \emph{bi-finite correspondence}) $\Gamma$ is defined to be an effective cycle in $X\times X$ with $\Q$-coefficients of dimension $d$  such that the projection from $\Gamma$ to the first (resp. each factor) is finite.
\end{defi}
The above notions are generalizations of the notions of endomorphisms and finite endomorphisms.
\begin{defi}
We call a cohomological correspondence $c\in \sC(X,X)$ \emph{finite} (resp. \emph{bi-finite}) if it takes form $$c=cl(\Gamma)$$ where $\Gamma$ is a finite (resp. bi-finite) correspondence.
\end{defi}

Here are some examples of bi-finite cohomological correspondences which I feel especially interesting:
\begin{exe}[Random product]
	Let $f_1,\dots, f_m$ be finite endomorphisms on $X$. Define $$c:=\sum_{i=1}^{m}a_if_i+\sum_{i=1}^{m}b_i{^{\top}f_i}$$
	with $a_i, b_i\in \Q_{\geq 0}.$ It is clear that $c$ is bi-finite.

	If $\sum_{i=1}^m(a_i+b_i)=1$, we may think that $c$ is a random product of $f_1,\dots, f_m$ and their transports $^{\top}f_1,\dots, ^{\top}f_m$ with probability $a_1,\dots, a_m, b_1,\dots, b_m.$ 
	One may rely the study of $c^n$ to the random product of matrices. 
\end{exe}

\begin{exe}[Extensions of endomorhisms]Let $f: X\to X$ be a surjective endomorphism. Let $Y$ be a smooth projecitve variety of pure dimension $d$.
	Let $\pi:Y\to X$ be a finite surjective morphism. Let $c_1,\dots, c_m$ be all irreducible components of $(\pi\times \pi)^{-1}(f).$ 
	For $a_i\in \Q_{\geq 0}, i=0,\dots, m$, the cohomological correspondence $c$ induced by $\sum_{i=1}^ma_ic_i$ is bi-fnite. 
	Assume that $$f\subseteq \Supp\,(\pi\times \pi)_*c.$$ Then by Lemma \ref{lemeffectvecom} and the projection formula, we may check that 
	$\beta_i(c)=\beta_i(f)$ for every $i=0,\dots, d.$
\end{exe}

For $c\in \sC(X,X)$.
Define the \emph{Cohomologcial polygon} for $c$ to be the minimal concave function 
$$CP_{c}: [0,2d]\to \R\cup\{-\infty\}$$ such that $$CP_{c}(j)\geq \log \alpha_j(c)$$ for all $j=0,\dots,2d.$

\medskip

Similarly, define the \emph{numerical polygon} for $c$ is the minimal concave function 
$$NP_c: [0,2d]\to \R\cup \{-\infty\}$$ such that $$NP_c(2i)\geq \log \beta_i(c)$$ for all $i=0,\dots,d.$

\medskip
As the numerical equivalence is weaker than the cohomological one, 
we have  $NP_c\leq CP_c.$
The following result shows that in the finite case, the above two polygons are indeed coincide.
\begin{thm}[=Theorem \ref{thmcpequnp}]\label{thmcpequnpintro}If $c$ is finite, then $NP_c= CP_c.$
\end{thm}

\subsection*{Acknowledgement}
The author would like to thank Xinyi Yuan and Qizheng Yin for helpful discussion. I especially thank Fei Hu and Tuyen Trung Truong for their helpful comments for the first version of the paper. I especially thank Weizhe Zheng, he simplified my original proof of Theorem \ref{thmcpequnp}.

\subsection{Notations}
\begin{points}
	\item[$\bullet$] Let $V$ be a finitely dimensional vector space over a field $K$ of characteristic $0$ and $f: V\to V$ be an endomorphism.
		We denote by $\Sp(f)$ the spectrum of $f$, i.e.  the multi-set of eigenvalues of $f$ and $P(f)$ the characteristic polynomial. 
			If $K$ is $\R$ or $\C$ and $V\neq \{0\}$, we denote by $$\rho(f):=\max\{|b||\,\, b\in \Sp(f)\}$$ the spectral radius of $f$.  For the convenience, if $V=\{0\}$, we define $\rho(f)=0.$
			If we fix any norm $\|\cdot\|$ and denote by $\|f^n\|$ the operator norm of $f^n$ for $n\geq 0$, we have $$\rho(f)=\lim_{n\to \infty} \|f^n\|^{1/n}.$$

\item[$\bullet$] A \emph{variety} is a separated reduced scheme of finite type over a field.
\item[$\bullet$] For any projective variety $Y$ of dimension $d$ and 
$i=0,\dots, d$, denote by $N_i(X)$ the group of numerical cycles of dimension $i$. 
Denote by $$N_i(X)_{\R}:=N_i(X)\otimes \R.$$
If $Y$ is further smooth, denote $$N^i(Y):=N_{d-i}(Y) \text{ and } N^i(Y)_\R:=N_{d-i}(Y)_\R.$$
	\end{points}

\section{Cohomological and numerical eigenvalues}
Let $X$ be a smooth projective variety over a field $\bk$ with pure dimension $\dim X=d.$
Let $l$ be a prime number  with $l\neq \Char\, \bk.$ We fix an embedding $\tau: \Q_l\hookrightarrow \C$ and view $\Q_l$ as a subfield of $\C$ via $\tau.$ Let $$G(X):=\{i=0,\dots, 2d|\,\, H^i(X_{\overline{\bk}},\Q_l)\neq 0\}.$$ By Poincar\'e duality, $i\in G(X)$ if and only if $2d-i\in G(X).$

\medskip

Denote by $\sC(X,X)_{\Z}$ the space of degree $0$ cohomological correspondences (for the $\Q_l$-cohomology) from $X$ to itself i.e. the image of the cycle class map  $$cl:\CH^d(X\times X)\to H^d(X_{\overline{\bk}}\times X_{\overline{\bk}})(d).$$  With the composition $\circ$, it forms a ring. Denote by $\sC(X,X):=\sC(X,X)_{\Z}\otimes \Q.$ For every endomorphism $f: X\to X$, we still denote by $f$ its graph in $X\times X$ and also the cohomological correspondences induced by it. 

\medskip
\subsection{Cohomological and numerical spectrum}
For $c\in \sC(X,X)$ and $j=0,\dots, 2d$, define
$$\Sp_j(c):=\Sp(c^*: H^j(X_{\overline{\bk}},\Q_l)\to H^j(X_{\overline{\bk}},\Q_l))$$ and 
$$P_j(c):=P(c^*: H^j(X_{\overline{\bk}},\Q_l)\to H^j(X_{\overline{\bk}},\Q_l)).$$
For $j\not\in G(X)$, we have $\Sp_j(c)=\emptyset$ and $P_j=1.$

\begin{fact}\label{factrationa}For every $j=0,\dots, 2d$,
	$P_j(c)$ has coefficients in $\Q$.  In particular, all elements of $\Sp_j(c)$ are algebraic.  
Moreover, if $c\in \sC(X,X)_{\Z}$, then $P_j(c)$ has coefficients in $\Z$. In particular, all elements of $\Sp_j(c)$ are algebraic integers. 
	\end{fact}
\proof
By the smooth-proper base change of \'etale cohomology and the spreading out argument, we may assume that $\bk=\F_q$ is a finite field.  
We conclude the proof by \cite[Theorem 2]{Katz1974}.
\endproof

For $c\in \sC(X,X)$, and $i=0,\dots, d$, define
$$\Sp_i^{num}(c):=\Sp(c^*: N^i(X)\otimes \R\to N^i(X)\otimes \R).$$
and 
$$P_j^{num}(c):=P(c^*: N^i(X)\otimes \R\to N^i(X)\otimes \R).$$
It is clear that $P_j^{num}(c)$ has coefficient in $\Q$. Moreover if $c\in \sC(X,X)_{\Z}$, then $P^{num}_i(c)$ has coefficients in $\Z$. 

\medskip

For $i=0,\dots, d$, let $A^{2i}(X)$ be the $\Q_l$-subspace of $H^{2i}(X_{\overline{\bk}},\Q_l)(i)$
generated by $cl(\CH^{i}(X)).$ It is invariant under $c^*$ and $c_*$ for every $c\in \sC(X,X).$
There is a natural morphism $$\phi:A^{2i}(X)\to N^i(X)\otimes \Q_l$$ as follows:
for every element $u\in A^{2i}(X)$, write $u=\sum_{t=1}^ma_t cl(Z_t)$ where $a_t\in \Q_l$ and $Z_t\in \CH^i(X).$
Define $\phi(u):=\sum_{t=1}^ma_t [Z_t]$ where $[Z_t]$ is the numerical class of $Z_t$.
It is well-defined. To see this, we only need to show that $\sum_{t=1}^ma_t cl(Z_t)=0$ implies that $\sum_{t=1}^ma_t [Z_t]=0.$
This is clear, as for every $W\in \CH^{d-i}(X)$, $$((\sum_{t=1}^ma_t [Z_t])\cdot W)=\langle (\sum_{t=1}^ma_t cl(Z_t))\cdot cl(W)\rangle=0.$$
It is clear the $\phi$ is surjective. Hence $N^i(X)\otimes \Q_l$ is a sub-quotient of $H^{2i}(X_{\overline{\bk}},\Q_l)(i).$ This implies that 

 $$\Sp_i^{num}(c)\subseteq \Sp_{2i}(c),\footnote{Here the inclusion is the inclusion for multi-sets i.e. for every element $a\in \Sp_i(c)$, it is in $\Sp_{2i}(c)$ with larger or equal multiplicity.}$$

\medskip

Easy to see the following fact.
\begin{fact}\label{factdecomcor}Let $X_1, X_2$ be smooth projective varieties of pure dimension $d$ and $X:=X_1\sqcup X_2.$
	Let $c_1\in \sC(X_1,X_1)$ and $c_2\in \sC(X_2,X_2)$ be bi-finite cohomological correspondences. 
	Then $c:=c_1\sqcup c_2\in \sC(X,X)$ is bi-finite. Moreover, we have $$\Sp_j(c)=\Sp_j(c_1)\sqcup \Sp_j(c_2)$$ for every $j=0,\dots, 2d$ and 
	$$\Sp_i^{num}(c)=\Sp_i^{num}(c_1)\sqcup \Sp_i^{num}(c_2)$$ for $i=0,\dots, d.$\footnote{Here the disjoint union is the disjoint union for multi-sets e.g. $\{a,a,b\}\sqcup \{a,c\}=\{a,a,a,b,c\}$.}
\end{fact}

If $f:X\to X$ is a surjective endomorphism of $X$, then there is $m\geq 1$ such that $f^m$ maps every irreducible component of $X$ to itself.
This is not true for general correspondence.
Using Fact \ref{factdecomcor}, many problems on $f$ can be reduced to the case where $X$ is irreducible.

\subsection{Cohomological spectral radius}\label{subsectioncohr}
Fix any embedding $\tau: \Q_l\hookrightarrow \C$ and view $\Q_l$ as a subfield of $\C$ via $\tau.$
For $c\in \sC(X,X)$ and $j=0,\dots, 2d$, 
define
 $$\alpha_j(c):= \rho(c^*: H^j(X_{\overline{\bk}},\Q_l)\otimes_{\Q_l}\C\to H^j(X_{\overline{\bk}},\Q_l)\otimes_{\Q_l}\C).$$ 
 By Fact \ref{factrationa},
 $\alpha_j(c)$ does not depend on the choice of $\tau.$
  We call all $\alpha_j(c)$ the $i$-th \emph{cohomological spectral radius} of $c$.
 Note that, if $j\not\in G(X)$, then  $\alpha_j(c)=0$.

\medskip

We say the $c$ has Condition (A) if there is a unique $i\in \{0,\dots, 2d\}$ such that 
$$\alpha_i(c)=\max_{j=0,\dots, 2d} \alpha_j(c).$$

Denote by $\Delta\in \sC(X,X)$ the diagonal of $X\times X.$
\begin{lem}\label{lemconadiag}
	If Condition (A) holds for $c\in \sC(X,X),$ then 
	$$\limsup_{n\to \infty} \langle c^n, \Delta\rangle^{1/n}=\max_{j=0,\dots,2d}\alpha_j(c).$$
	\end{lem}
\proof 
We have 
\begin{equation}\label{equtracinter}\langle c^n, \Delta\rangle=\sum_{j=0}^{2d}(-1)^j\Tr((c^n)^*: H^j(X_{\overline{\bk}},\Q_l)\otimes_{\Q_l}\C).
	\end{equation}
It is clear that for $j=0,\dots, d.$
\begin{equation}\label{equtrlimsup}
	\limsup_{n\to \infty}|\Tr((c^n)^*: H^j(X_{\overline{\bk}},\Q_l)\otimes_{\Q_l}\C)|^{1/n}=\alpha_j(c).
\end{equation}
There is a unique $i\in \{0,\dots, 2d\}$ such that 
\begin{equation}\label{equtrconda}
\alpha_i(c)=\max_{j=0,\dots, 2d} \alpha_j(c).
\end{equation}
Combine (\ref{equtrlimsup}) with (\ref{equtrconda}), we get 
$$\limsup_{n\to \infty}|\Tr((c^n)^*: H^i(X_{\overline{\bk}},\Q_l)\otimes_{\Q_l}\C)^{1/n}|=\alpha_i(c)$$
and there is $\alpha<\alpha_i(c)$ such that for $n\gg 0$, we have 
$$|\sum_{j\neq i}\Tr((c^n)^*: H^j(X_{\overline{\bk}},\Q_l)\otimes_{\Q_l}\C)|\leq \alpha^n.$$
So we get 
$$\limsup_{n\to \infty}(\sum_{j=0}^{2d}\Tr((c^n)^*: H^j(X_{\overline{\bk}},\Q_l)\otimes_{\Q_l}\C))^{1/n}=\alpha_i(c),$$
which concludes he proof.
\endproof

\subsection{Finite cohomological correspondence}
\begin{defi}We call a cohomological correspondence $c\in \sC(X,X)$ \emph{effective} if it takes form $$c=cl(\Gamma)$$ where $\Gamma$ is an effective cycle in $X\times X$ of dimension $d.$ 
\end{defi}
If $c\in \sC(X,X)$ is effective, then its transport $^{\top}c$ is also effective.
If $c_1,c_2\in \sC(X,X)$ are effective and $r_1,r_2\in \Q_{\geq 0}$, $r_1c_1+r_2c_2$ is effective.

\begin{defi}We call a cohomological correspondence $c\in \sC(X,X)$ finite (resp.\emph{bi-finite}) if it takes form $$c=cl(\Gamma)$$ where $\Gamma$ is an effective cycle in $X\times X$ of dimension $d$ such that the projection from $\Gamma$ to the first (resp. each) factor is finite.
	\end{defi}
It is clear that bi-finite correspondences are effective.
We have the following basic properties:
\begin{points}
	\item If $c\in \sC(X,X)$ is bi-finite, then its transport $^{\top}c$ is also bi-finite.
\item If $c_1,c_2\in \sC(X,X)$ are finite (resp. bi-finite), then $c_1\circ c_2$ is finite (resp. bi-finite).
\item If $c_1,c_2\in \sC(X,X)$ are finite (resp. bi-finite) and $r_1,r_2\in \Q_{\geq 0}$ then $r_1c_1+r_2c_2$ is finite (resp. bi-finite).
\end{points}
In particular, for every if $c\in \sC(X,X)$ is finite (resp. bi-finite), then $c^n$ is finite (resp. bi-finite) for every $n\geq 0.$

\medskip

For a bi-finite $c\in \sC(X,X)$, we say the $c$ has a \emph{bi-finite inverse} if it has an inverse $c^{-1}\in \sC(X,X)$ which is also bi-finite.
If a bi-finite $c\in \sC(X,X)$ has a bi-finite inverse, then for every $n\in \Z$, $c^n$ is bi-finite.

\begin{exe}\label{exeendo}
	Let $f: X\to X$ be a surjective endomorphism. By \cite[Lemma 5.6]{fa}, $f$ is finite.
	If we view $f$ as a cohomological correspondence in $C(X,X)$, then $f$ is bi-finite.
	As $$f^{-1}=(\deg f)^{-1}{^{\top}f}\in C(X,X),$$
	$f^{-1}$ is bi-finite.
\end{exe}

\begin{lem}\label{lembifinitecone}Let $c_1,c_2,e\in \sC(X,X). $ Assume that $c_1,c_2$ are bi-fnite and $e$ is effective. Then $c_1\circ e\circ c_2$ is effective. 
	\end{lem}
\proof
It is clear that for every bi-finite $c\in \sC(X,X)$ and effective $e\in \sC(X,X)$, $e\circ c$ is effective. 
Taking transport, we get that $c\circ e$ is also effective. Then we conclude the proof.
\endproof


\medskip

\subsection{Numerical spectral radius}\label{subsectionnumr}
For $c\in \sC(X,X)$, define $$\beta_i(c):=\rho(c^*: N^i(X)_{\R}\to N^i(X)_\R),$$
and call it the $i$-th \emph{numerical spectral} of $c$.
As the cohomological equivalence is finner than the numerical equivalence,
we have 
\begin{equation}\label{equnumlch}
	\beta_i(c)\leq \alpha_{2i}(c).
\end{equation}

\medskip

Let $\Psef^i(X)$ be the pseudo-effective cone in $N^i(X)$.
By \cite[Corollary 3.3.7]{Dang2020}, $\Psef^i(X)$
is a closed convex cone with non-empty interior and moreover it is salient i.e. $\Psef^i(X)\setminus \{0\}$ is convex. Let $L$ be an ample line bundle on $X$. 
By \cite[Proposition 3.3.6]{Dang2020} and \cite[Theorem 3.3.3 (ii),(v)]{Dang2020}, $L^i$ is contained in the interior of $\Psef^i(X)$ and moreover for every $u\in \Psef^i(X)$, $(u\cdot L^{d-i})=0$ if and only if $u=0.$
We have a norm $\|\cdot\|_L$ on $N^i(X)_{\R}$ as follows:
for $u\in N^i(X)_{\R}$, $$\|u\|_L:=\inf\{(u^+\cdot L^{d-i})+(u^-\cdot L^{d-i})|\,\, u^+,u^-\in \Psef^i(X), u=u^+-u^-\}.$$
Easy to check that $\|\cdot\|_{L}$ is a norm and for every $u\in \Psef^i(X)$, $\|u\|_L=(u\cdot L^{d-i}).$
If $c\in \sC(X,X)$ is finite (resp. bi-finite), then $c_*$ (resp. both $c_*$ and $c^*$) preserves $\Psef^i(X)$.
The following lemma give another description of $\beta_i(c).$
\begin{lem}\label{lemdyndegint}If $c\in \sC(X,X)$ is finite, then for every $i=0,\dots, d$, we have
	$$\beta_i(c)=\lim_{n\to \infty}((c^n)^*(L^i)\cdot L^{d-i})^{1/n}.$$
\end{lem}
\proof
By duality, we have $$\beta_i(c)=\rho(c_*:N^{d-i}(X)\otimes\R\to N^{d-i}(X)\otimes \R).$$

For every $u\in \Psef^{d-i}(X)$, there is $C>0$ such that $CL^{d-i}-u\in \Psef^i(X).$
Then for every $n\geq 0$, we have $(c^n)_*(CL^{d-i}-u)\in \Psef^{d-i}(X).$ So we have
$$(C\|(c^n)_*L^{d-i}\|_L)^{1/n}\geq \|(c^n)_*u\|_L^{1/n}$$ for every $n\geq 0.$
As $\Psef^{d-i}(X)$ has non-empty interior, the above inequality shows that 
$$\lim_{n\to \infty}((c^n)^*(L^i)\cdot L^{d-i})^{1/n}=\lim_{n\to \infty}\|(c^n)_*L^{d-i}\|_L^{1/n}$$
$$=\rho(c^*: N^{d-i}(X)_{\R}\to N^{d-i}(X)_\R)=\beta_i(c).$$
\endproof

Let $\pi_i:X\times X\to X, i=1,2$ be the projection to the first and the second coordinate. Set $L_i:=\pi_i^*L.$
Then $L_1+L_2$ is ample on $X\times X.$

\begin{lem}\label{lemeffectvecom}Assume that $c_1,c_2\in \sC(X,X)$ are bi-finite and $c_1-c_2$ is effective. Then for every $i=0,\dots, d$, we have $\beta_i(c_1)\geq \beta_i(c_2).$
\end{lem}

\proof
Set $e:=c_1-c_2.$
Then for $n\geq 0$, $$c_1^n-c_2^n=\sum_{l=0}^{n-1}c_2^l\circ e\circ c_1^{n-1-l}$$ is effective.
So we have $((c_1^n-c_2^n)\cdot L_1^i\cdot L_2^{d-i})\geq 0.$
We conclude the proof by Lemma \ref{lemdyndegint}.
\endproof

Denote by $\Delta\in \sC(X,X)$ the diagonal of $X\times X.$

\begin{lem}\label{leminterdegree}Assume that $c\in \sC(X,X)$ is finite. Then we have 
$$\limsup_{n\to \infty}\langle c^n, \Delta \rangle^{1/n}\leq \max_{i=0}^d\beta_i(c).$$
	\end{lem}
\proof
Consider the linear function $I: N^d(X\times X)_{\R}\to \R$ defined by $u\mapsto (u\cdot \Delta).$
There is $C>0$ such that $$I(u)\leq C\|u\|_{L_1+L_2}$$ for every $u\in N^d(X\times X)_{\R}.$

We still denote by $c^n$ the numerical class in $N^d(X\times X)$ induced by $c^n$.
Then we have 
\begin{equation}\label{equcnbonorm}
	\langle c^n, \Delta \rangle=(c^n\cdot \Delta)=I(c^n)\leq C\|c^n\|_{L_1+L_2}.
\end{equation}
As $$\|c^n\|_{L_1+L_2}=\sum_{i=0}^d \binom{d}{i}((c^n)^*{L^i}\cdot L^{d-i}),$$
we get $$\lim_{n\to \infty}\|c^n\|_{L_1+L_2}^{1/n}=\max_{i=1}^d\lim_{n\to \infty}((c^n)^*{L^i}\cdot L^{d-i})^{1/n}=\max_{i=0}^d\beta_i(c).$$
We conclude the proof by (\ref{equcnbonorm}).
\endproof

By Lemma \ref{leminterdegree}, Lemma \ref{lemconadiag} and (\ref{equnumlch}), we get the following consequence.
\begin{cor}\label{corassuaent}Assume that $c\in \sC(X,X)$ is finite and has Condition (A). Then we have
	$$\max_{j=0}^{2d}\alpha_i(c)=\max_{i=0}^d\beta_i(c).$$
\end{cor}

Indeed, this corollary is true without condition (A). By smooth-proper base change, we only need to treat the case where $\bk$ is a finite field. On the other hand, when $\bk$ is a finite field, this corollary was proved by Truong in \cite[Theorem 2.1]{Truong2016} without Condition (A). However, Truong's proof strongly relies on Deligne's theorem for Weil's  Riemann Hypothesis \cite{Deligne1974}. Our proof is much elementary and for our purpose, it's not hard to get condition (A), so we follow this more elementary way.

\subsection{Cohomologcial and  numerical polygons}
Let $c\in \sC(X,X)$.
Define the \emph{Cohomologcial polygon} for $c$ to be the minimal concave function 
$$CP_{c}: [0,2d]\to \R\cup\{-\infty\}$$ such that $$CP_{c}(j)\geq \log \alpha_j(c)$$ for all $j=0,\dots,2d.$
Observe that $c$ has condition (A) if and only if there is a unique $j\in \{0,\dots, 2d\}$ such that $CP_{c}(j)$ takes the maximal value. 

\medskip

Similarly, define the \emph{numerical polygon} for $c$ is the minimal concave function 
	$$NP_c: [0,2d]\to \R\cup \{-\infty\}$$ such that $$NP_c(2i)\geq \log \beta_i(c)$$ for all $i=0,\dots,d.$

\medskip

We note that $$\max_{i=0,\dots d} NP_c(2i)=\max_{i=0,\dots,d}\beta_i(c)$$
and $$\max_{j=0,\dots 2d} CP_c(j)=\max_{j=0,\dots,2d}\alpha_j(c).$$
By (\ref{equnumlch}), we have 
\begin{equation}\label{equcpnp}
	NP_c\leq CP_c.
	\end{equation}

\medskip
By Lemma \ref{lemdyndegint}, if $c$ is finite, then we have $CP_c\geq NP_c$ on $[0, 2d].$
Moreover, if $c$ is bi-finite, then $NP_c>-\infty.$
\begin{thm}\label{thmcpequnp}If $c$ is finite, then $NP_c= CP_c.$
	\end{thm}
	
\proof
By the smooth-proper base change of \'etale cohomology and the spreading our argument, we may assume that $\bk=\F_q$ is a finite field.

We may approximate $c$ by $c_m:=c+m^{-1}\Delta$ with $m\to \infty.$ 
For each $m\geq 1$, as $c=c_m-m^{-1}\Delta$ is effective. As $m^{-1}\Delta$ is bi-finite, and $m^{-1}\Delta$
commutes with $c_m$, the proof of 
Lemma \ref{lemeffectvecom} indeed shows that $$\beta_i(c_m)\geq \beta_i(m^{-1}\Delta)>0$$ for every $i=0,\dots, d$.
After replacing $c$ by $c_m, m\geq 1$, we may assume that $\beta_i(c)>0$ for every $i=0,\dots, d.$




\medskip

Denote by $\Phi_q$ the cohomological correspondence in $\sC(X,X)$ induced by the $q$-Frobenius on $X$. So $\Phi_q$ is bi-finite and $\Phi_q^{-1}=q^{-d}{^{\top}\Phi_q}.$ 
Deligne's theorem \cite{Deligne1974} implies that all eigenvalues of $\Phi_q^*: H^{j}(X_{\overline{\bk}},\Q_l)\otimes \C\to H^{j}(X_{\overline{\bk}},\Q_l)\otimes \C$ has norm $q^{j/2}.$
For $s\in \Z, t\in \Z_{\geq 0}$, define 
$c_{s,t}:=\Phi_q^s\circ c^t.$ As $\Phi_q$ commutes with $c,$ we have the following properties: for $x\in [0,2d]$,
\begin{points}
	\item $NP_{c_{s,t}}(x)=(\frac{1}{2}s\log q)x+tNP_c(x);$
	\item $CP_{c_{s,t}}(x)=(\frac{1}{2}s\log q)x+tCP_c(x).$
\end{points}

\medskip
These $NP_{c_{s,t}}$, $CP_{c_{s,t}}$ are continuous, piece-wise linear and concave functions on $[0,2d]$ satisfying 
$$NP_{c_{s,t}}\leq CP_{c_{s,t}}.$$
 
\medskip

For any continuous, piece-wise linear and concave function $h:[0,2d]\to \R$. For $x\in [0,2d]$, denote by $d^-h(x),d^+h(x)$ the left and the right derivation of $h.$ For the convenience, we set $d^-h(0)=\infty$ and $d^+h(2d)=-\infty.$
Denote by $V(h)$ the set of the points $x\in [0,2d]$ satisfying $d^-h(x)>d^+h(x)$.
Then $V(h)$ is a finite set containing $\{0,2d\}.$ It is clear that for every $s\in \Z, t\in \Z_{\geq 0}$,
$V(CP_{c_{s,t}})=V(c)\subseteq \{0,\dots, 2d\}.$ 

\medskip

We only need to prove that   
$$NP_{c_{s,t}}\geq CP_{c_{s,t}}.$$
For this, we only need to show that for every 
$x\in V(c)$, $CP_c(x)\leq NP_c(x).$
Pick $(s,t)\in \Z\times \Z_{>0}$ such that $$s/t\in (-2d^-CP_c(x)/\log q, -2d^+CP_c(x)/\log q).$$
We have
$$d^-CP_{c_{s,t}}(x)>0 \text{ and } d^+CP_{c_{s,t}}(x)<0.$$
Hence $CP_{c_{s,t}}$ takes maximal value only at $x.$ 
Then $c_{s,t}$ has  condition (A).
By Corollary \ref{corassuaent},
we get $$CP_{c_{s,t}}(x)=\max CP_{c_{s,t}}=\max NP_{c_{s,t}}.$$
We claim that $NP_{c_{s,t}}$ also takes maximal value only at $x$. Otherwise, there is $y\in [0,2d]\setminus \{x\}$ such that 
$$NP_{c_{s,t}}(y)=\max NP_{c_{s,t}}=CP_{c_{s,t}}(x)>CP_{c_{s,t}}(y)\geq NP_{c_{s,t}}(y),$$ which is a contradiction. 
Then we get $$CP_{c_{s,t}}(x)=NP_{c_{s,t}}(x).$$
By (i) and (ii) above, we get 
$$CP_{c}(x)=NP_{c}(x),$$
which concludes the proof.
	\endproof

\subsection{Case of endomorphisms}
In this section, we assume that $X$ is {\bf irreducible}.
Let $f: X\to X$ be a surjective endomorphism. 
We still denote by $f$ the cohomological correspondence in $\sC(X,X)$ induced by $f$.

\medskip

For any bi-finite $c\in \sC(X,X),$ we say that $c$ is \emph{numerically log-concave}, if 
\begin{equation}\label{equationcorlogconcave}
	NP_c(2i)=\log \beta_i(c).
\end{equation}

By \cite[Theorem 1.1(3)]{Truong2020} (see also \cite{Dinh2005,Dang2020}), if $X$ is irreducible, the sequence $\beta_i(f), i=0,\dots, d$ is log-concave.
So $f$ is numerically log-concave. Indeed, as the transport of the graph of $f^n$ is irreducible for all $n\geq 0$,
\cite[Theorem 1.1(3)]{Truong2020} also implies the $^{\top}f$ is numerically log-concave.

\medskip

The following gives an example of bi-finite cohomological correspondence $c\in \sC(X,X)$ which is not numerically log-concave.
\begin{exe}
Let $f_i, i=1,2$ be two surjective endomorphisms on $\P^d$ of algebraic degree $q_i\geq 1$ i.e. $f_i^*O(1)=O(q_i)$.
Define $c:=\frac{1}{2}(f_1+f_2)\in \sC(X,X).$ We may think $c$ as a random product of $f_1, f_2$ independently of probability $1/2$ for each $i.$
Set $L:=O(1).$ For $i=0,\dots,d$, we have 
\begin{equation}\label{equationcompbetac}\beta_i(c)=\lim_{n\to\infty}((c^n)^*L^i\cdot L^{d-i})^{1/n}=\lim_{n\to\infty}(2^{-n}(q_1^i+q_2^i)^n)^{1/n}=(q_1^i+q_2^i)/2.
	\end{equation}
Cauchy inequality implies that for every $x>0$, $i\geq 1$, we have $$(1+x^{i-1})(1+x^{i+1})\geq (1+x^i)^2,$$
and the equality holds if and only if $x=1.$
Applying the above inequality, we get that for every $i=1,\dots, d-1$, we get
$$\beta_{i-1}(c)\beta_{i+1}(c)\geq \beta_i(c)^2,$$
and  the  equality holds if and only if $q_1=q_2.$
In particular, if $q_1\neq q_2$, then 
	$$NP_c(2i)>\log \beta_i(c)$$
	for all $i=1,\dots, d-1.$
	\end{exe}
	
\medskip

In general the sequence $\alpha_i(f), i=0,\dots, 2d$ may not be log-concave even when $X$ is irreducible. For example, it is possible that $H^i(X_{\overline{\bk}},\Q_l)$ vanishes for some odd $i$. Then this case $\alpha_i(f)=0$. As shown in the following example, even when $H^i(X_{\overline{\bk}},\Q_l)\neq 0$ for every $i=0,\dots, 2d$, $\alpha_i(f), i=0,\dots, 2d$ could not be log-concave.
\begin{exe}
Let $X=\P^1\times E$ where $E$ is an elliptic curve. Let $g: \P^1\to \P^1$ be the square map $[x:y]\mapsto [x^2:y^2].$
Set $f:=g\times \id: X\to X.$ 
Easy to see that $H^i(X_{\overline{\bk}},\Q_l)\neq 0$ for all $i=0,\dots, 4$. Easy to compute that $\alpha_0(f)=1, \alpha_1(f)=1, \alpha_2(f)=2$.
As $$1=\alpha_1(f)^2<\alpha_0(f)\alpha_1(f)=2,$$ $\alpha_i(f), i=0,\dots, 4$ is not log-concave.
\end{exe}

\medskip

For a finite multi-set $A\subseteq \C^*$ and $b\in \C^*$, write $$A^{-1}=\{a^{-1}|\,\,a\in A\}, bA:=\{ba|\,\, a\in A\}
\text{ and } |A|=\{|a||\,\, a\in A\}.$$

Now we fix an embedding $\tau: \Q_l\hookrightarrow \C.$
By Poincar\'e duality, we have the following facts:
For every $j\in G(X), i=0,\dots, d$, we have 
	\begin{points}
\item $0\not\in \Sp_j(f)$ and $0\not\in \Sp(f^*: N^i(X)_{\R}\to N^i(X)_{\R})$,
\item $\Sp_j(^{\top}f)=(\deg f)\Sp_j(f)^{-1}$ and $$\Sp((^{\top}f)^*: N^i(X)_{\R}\to N^i(X)_{\R})=(\deg f)\Sp(f^*: N^i(X)_{\R}\to N^i(X)_{\R})^{-1};$$
\item $\Sp_j(f)=\Sp_{2d-j}(^{\top}f)$ and $$\Sp(f^*: N^i(X)_{\R}\to N^i(X)_{\R})=\Sp((^{\top}f)^*: N^{d-i}(X)_{\R}\to N^{d-i}(X)_{\R}).$$
\end{points}
For $j\in G(X)$, set $$\alpha^-_j(f):=\min|\Sp_j(f)|.$$
For $i=0,\dots, d$, set $$\beta^-_i(f):=\min|\Sp(f^*: N^i(X)_{\R}\to N^i(X)_{\R})|.$$

By Fact \ref{factrationa}, $\rho^-_j(f)$ does not depend on the choice of the embedding $\tau$.
By the above facts (ii) and (iii), we get the following result.
\begin{pro}\label{prospinf}
	For every $j\in G(X)$, we have 
	$$\alpha^-_j(f)=\frac{\deg(f)}{\alpha_{2d-j}(f)} \text{ and } \beta^-_i(f)=\frac{\deg(f)}{\beta_{d-i}(f)}.$$
	\end{pro}

By Theorem \ref{thmcpequnp}, we get the  following consequence.
\begin{cor}\label{corendosum}For every $j\in G(X)$, we have 
$$NP_f(2d)-NP_f(2d-j)\leq \log \alpha^-_j(f)\leq \log\alpha_j(f)\leq NP_f(j)$$
Moreover, for every $i=0,\dots, d$, we have 
$$NP_f(2d)-NP_f(2d-2i)=\log \alpha^-_{2i}(f)=\log\beta^-_i(f)$$
and 
$$\log\beta_i(f)=\log \alpha_{2i}(f)=NP_f(2i).$$
Note that $\log\deg f=NP_f(2d).$
\end{cor}
\medskip

In the end, we prove the ``moving target" version of  Corollary \ref{corfixpoint}.
\begin{pro}\label{promovingtaget}Assume that $f: X\to X$ is an int-amplified amplified endmorphism. Let $L$ be an ample line bundle on $X$.
	Let $h_n: X\to X, n\geq 0$ be a sequence of endomorphisms of $X$ with 
	$$\limsup_{n\to \infty}(h_n^*L\cdot L^{d-1})^{1/n}<\beta_{d}(f)/\beta_{d-1}(f).$$
	Then for $n\gg 0$, $f^n$ intersects $h_n$ properly in $X\times X$. Moreover, for every $\epsilon>0$ we have 
	$$\#\{f^n(x)=h_n(x)\}=q^{dn}+o((\beta_d\beta_{d-1})^{(1+\epsilon)n})$$
	counting with multiplicity.
	\end{pro}
\proof
Let $\pi_i: X\times X\to X, i=1,2$ be the $i$-th projection. 
Set $\mu_d:=\beta_{d}(f)/\beta_{d-1}(f).$ As $f$ is int-amplified amplified, $\mu_d>1.$
By Corollary \ref{corendosum}, every eigenvalue of $f^*: N^1(X)\otimes \R\to N^1(X)\otimes \R$ has norm at least $\mu_d.$
Hence every eigenvalue of $f_*: N^{d-i}(X)\otimes \R\to N^{d-i}(X)\otimes \R$ has norm at least $\mu_d.$
So there is $\delta>0,$ such that for every $Z\in \Psef^{d-i}(X)$ and $n\geq 0$ we have 
\begin{equation}\label{equationnonegrow}(f_*^n(Z)\cdot L)=\|f_*^n(Z)\|_L\geq \delta\mu_d^n\|Z\|_L.
	\end{equation}
There is a constant $A>0$, such that for for every $M\in N^1(X),Z\in N^{d-1}(X)$, we have
\begin{equation}\label{equationcomparenorm}(Z\cdot M)\leq A\|M\|_L\|Z\|_L.
\end{equation}
Pick $\eta\in (\limsup_{n\to \infty}(h_n^*L\cdot L^{d-1}), \mu_d).$
There is $B>0$ such that 
\begin{equation}\label{equationhndegup}
	(h_n^*L\cdot L^{d-1})<B\eta^n.
	\end{equation}

\medskip

We first prove that $f^n$ intersects $h_n$  properly for $n\gg 0$.
Assume that there is an irreducible curve $C'\subseteq f^n\cap h_n\subseteq X\times X$ for some $n\geq 0.$
Set $C:=\pi_1(C').$ We have $(f^n)^*L|_C=(h_n)^*L|_C.$ So we get 
$$(L\cdot f^n_*(C))=(h_n^*L\cdot C).$$
By (\ref{equationnonegrow}), (\ref{equationcomparenorm}) and (\ref{equationhndegup}), we get 
$$\delta\mu_d^n\|C\|_L\leq\|f^n_*(C)\|_L=(L\cdot f^n_*(C))=(h_n^*L\cdot C)\leq A(h_n^*L\cdot L^{d-1})\|C\|_L\leq AB\eta^n\|C\|_L.$$
Hence $$n\leq \frac{\log A+\log B-\log \delta}{\log\mu_d-\log \eta}.$$

\medskip

Fix a field embedding $\tau: \Q_l\hookrightarrow \C$.
For every $i=0,\dots, 2d$, $V_i:=H^i(X_{\overline{\bk}},\Q_l)\otimes_{\Q_l}\C.$
Fix a norm $\|\cdot\|$ on each $V_i$. For every endomorphism $g: V_i\to V_i$, denote by $\|g\|$ the operator norm of $g$.
There is a constant $D>0$ such that for every $i=0,\dots, 2d$ and $g: V_i\to V_i$, we have 
$$|\Tr(g)|\leq D\|g\|.$$
For $n\gg 0$, we have
$$\#\{f^n(x)=h(x)\}=(f^n\cdot h_n)=\sum_{i=0}^{2d}(-1)^i\Tr((^{\top}h)^*\cdot (f^n)^*: V_i\to V_i)$$
$$=\mu_d^n+\sum_{i=0}^{2d-1}(-1)^i\Tr((^{\top}h)^*\cdot (f^n)^*: V_i\to V_i).$$

We only need to bound $\sum_{i=0}^{2d-1}(-1)^i\Tr((^{\top}h)^*\cdot (f^n)^*: V_i\to V_i).$
For every $i=0,\dots, 2d-1$, we have 
$$|\Tr((^{\top}h)^*\cdot (f^n)^*|_{V_i})|\leq D\|(^{\top}h)^*|_{V_i}\|\|(f^n)^*|_{V_i}\|\leq D\|(^{\top}h)^*|_{V_i}\|\alpha_i(f)^{(1+\epsilon/2)n}.$$
As $\beta_i, i=1,\dots, d$ is increasing,
by Corollary \ref{corendosum}, for every $i=0,\dots, 2d-1$, we get $\alpha_i(f)\leq (\beta_{d-1}\beta_{d})^{1/2}.$ This concludes the proof.
\endproof

	\bibliography{dd}
\end{document}